\documentclass[10pt,preprint]{article}

\usepackage{amsmath, latexsym, amsfonts, amssymb, amsthm}
\usepackage{graphicx, color,hyperref,dsfont,epsfig,caption,wrapfig,subfig}
\usepackage{pstricks}

\usepackage{movie15}
\hypersetup{
    linktoc=page,
    linkcolor=red,          
    citecolor=blue,        
    filecolor=blue,      
    urlcolor=cyan,
    colorlinks=true           
}

\numberwithin{equation}{section}

\newtheorem{theorem}{Theorem}[section]
\newtheorem{lemma}[theorem]{Lemma}

\newtheorem{rem}[theorem]{Remark}
\newtheorem{remark}[theorem]{Remark}

\theoremstyle{definition}

\renewcommand{\tilde}{\widetilde}          
\DeclareMathSymbol{\leqslant}{\mathalpha}{AMSa}{"36} 
\DeclareMathSymbol{\geqslant}{\mathalpha}{AMSa}{"3E} 
\DeclareMathSymbol{\eset}{\mathalpha}{AMSb}{"3F}     
\renewcommand{\leq}{\;\leqslant\;}                   
\renewcommand{\geq}{\;\geqslant\;}                   

\newcommand{\R}{\mathbb{R}}

\newcommand{\N}{\mathbb{N}}

\newcommand{\E}{\mathds{E}}
\newcommand{\Pb}{\mathds{P}}
\newcommand{\ind}{\mathds{1}}

\usepackage{xparse}

\def\Bvec{\boldsymbol B}

\DeclareDocumentCommand \Pmp { m m o} {
\IfNoValueTF{#3}
{P_{#1}^{#2}}
{P_{#1}^{#2}\left(#3\right)}
}

\DeclareDocumentCommand \Emp { m m o} {
\IfNoValueTF{#3}
{E_{#1}^{#2}}
{E_{#1}^{#2}\left[#3\right]}
}

\DeclareDocumentCommand \Pbr { m m m m o } {
\IfNoValueTF{#5}
{P_{#1}^{#2\stackrel{#4}{\rightarrow} #3}}
{P_{#1}^{#2\stackrel{#4}{\rightarrow} #3}\left(#5\right)}
}

\DeclareDocumentCommand \Ebr { m m m m o } {
\IfNoValueTF{#5}
{E_{#1}^{#2\stackrel{#4}{\rightarrow} #3}}
{E_{#1}^{#2\stackrel{#4}{\rightarrow} #3}\left[#5\right]}
}

\def\eps{\varepsilon}

\def\bi{\begin{itemize}}
\def\ei{\end{itemize}}
\def\ni{\noindent}
\def\bnum{\begin{enumerate}}
\def\enum{\end{enumerate}}
\def\LB{\mathcal{B}} 
\def\<#1{\langle #1 \rangle}

\def\M{\mathbf{M}}

\def\p{\mathbf{p}}

\def\cF{\mathcal{F}}

\title{KPZ formula derived from Liouville heat kernel}

\author{ Nathana\"el Berestycki\footnote{University of Cambridge, Statistical Laboratory. Research supported in part by EPSRC grants EP/GO55068/1 and EP/I03372X/1.}, Christophe Garban\footnote{ENS, Lyon, CNRS, France. Partially supported by ANR grant MAC2 10-BLAN-0123.}, R\'emi Rhodes \footnote{Universit\'e Paris Est-Marne la Vall\'ee, LAMA and CNRS UMR 8050, France. Partially supported by grant ANR-11-JCJC  CHAMU.} \footnotetext[2]{Partially supported by grant ANR-11-JCJC  CHAMU.},
 Vincent Vargas \footnote{ENS Ulm, DMA, 45 rue d'Ulm,  75005 Paris, France. Partially supported by grant ANR-11-JCJC  CHAMU.} }

\date{\vspace{-5ex}}

\begin{document}

\maketitle

 
\begin{abstract}
In this paper, we establish the Knizhnik--Polyakov--Zamolodchikov (KPZ) formula of Liouville quantum gravity, using the heat kernel of Liouville Brownian motion. This derivation of the KPZ formula was first suggested by F. David and M. Bauer in order to get a geometrically more intrinsic way of measuring the dimension of sets in Liouville quantum gravity. We also provide a careful study of the (no)-doubling behaviour of the Liouville measures in the appendix, which is of independent interest. 
\end{abstract}
\footnotesize



\noindent{\bf Key words or phrases:}  KPZ formula, heat kernel, Liouville quantum gravity, Liouville measures, Gaussian multiplicative chaos.


\normalsize



\section{Introduction}
Let $\gamma \ge 0$ and consider a free field $X$ in $\R^2$   and the (formal) Riemannian metric tensor
\begin{equation}\label{tensor}
e^{\gamma X(x)}dx^2.
\end{equation}
The tensor \eqref{tensor} gives rise to a random geometry known in physics as (critical) Liouville quantum gravity (LQG for short); see \cite{cf:Da,DistKa,cf:DuSh,bourbaki,cf:KPZ,Nak,Pol,review} for a series of works both within the physics and mathematics literature on the subject. 
While a rigorous construction of the metric associated to \eqref{tensor} is still an open problem, it is believed that \eqref{tensor} provides the scaling limit of certain random planar maps coupled with statistical physics models, and then suitably embedded in the plane (for instance via circle packing). The parameter $\gamma$ in \eqref{tensor} is related to the model of statistical physics under consideration. See for example the survey \cite{bourbaki}. 

One of the striking features of Liouville quantum gravity is the Knizhnik--Polyakov--Zamolodchikov (KPZ) formula. This is a far-reaching identity, relating the `size' (dimension) of a given set $A \subset \R^2$ from the point of view of standard Euclidean geometry, to its counterpart from the point of view of the geometry induced by \eqref{tensor}. 
The original formulation of the KPZ formula was made in \cite{cf:KPZ} in the context of the light-cone gauge, see also \cite{cf:Da,DistKa} for a derivation in the conformal gauge.

Recently, a rigorous version of the KPZ formula has been established in \cite{cf:DuSh,Rnew10,Rnew4,Rnew12,Aru}   (see also \cite{benj} in the context of one dimensional Mandelbrot's cascades), which we recall now. Usually, on a metric space, one defines the Hausdorff dimension of a given set $A$ by studying the diameters of open sets required  to cover this set. Yet, as already mentioned above, in LQG the construction of such a metric space is still an open question: the existence of a distance associated to the metric tensor \eqref{tensor} has not been verified yet. Let us mention here  the  work \cite{MS} for some recent advances on this topic  in the case $\gamma=\sqrt{8/3}$. 
At the time of writing, well-defined objects associated to this metric tensor are the volume form $M_\gamma$, which may be
 defined by the theory of Gaussian multiplicative chaos \cite{cf:Kah}; and the Liouville Brownian motion, which is the natural diffusion process associated to \eqref{tensor}. This may be constructed via potential theory as in \cite{GRV1} (see  also \cite{berest} for the construction starting from one point via convolution techniques). 
 
 To go around the difficulty of formulating the notion of Hausdorff dimension in the absence of a well-defined distance, different notions of Hausdorff dimensions have been suggested. For instance, a measure-based notion was formulated in \cite{Rnew10,Rnew4}. Instead of working with the diameter of open sets covering the set $A$, one works with the measure of small Euclidean balls covering this set (see subsection \ref{backKPZ} for more details). Another formulation is suggested in \cite{cf:DuSh}: the authors work with a notion of expected box counting dimension, once again by measuring the size of Euclidean balls with the measure $M_\gamma$. Finally, yet another formulation is suggested in \cite{Aru} based on the Minkowski dimension.  In these three cases, this yields a notion of quantum dimension, for which  the papers \cite{cf:DuSh,Rnew10,Aru} establish the KPZ formula  relating the Euclidean Hausdorff dimension of a set $A$ to its quantum dimension. However, in both cases, as observed by \cite{David-KPZ}, this formulation of quantum dimension is not geometrically intrinsic to LQG: the definition of quantum dimension deeply relies on the Euclidean structure, since 
to compute the quantum dimension, one covers the set with Euclidean balls. 

A conceptually more satisfying, and geometrically more intrinsic, derivation of the KPZ formula was proposed in \cite{David-KPZ}, based on the heat kernel of Liouville Brownian motion. 
Roughly speaking, the idea is to measure how the natural diffusion associated to a metric ``fills" the fractal sets. More precisely consider a fractal set $A$ with a non trivial $2-2q$ Euclidean Hausdorff measure $\mathcal{H}^{2-2q}$, i.e. $0< \mathcal{H}^{2-2q}(A)  <\infty$.  Consider the standard heat kernel $p_t(x,y)=(2\pi t)^{-1} e^{-|x-y|^2/2t}$ of the planar Brownian motion. Then a straightforward computation shows that  
\begin{equation}\label{scalingeuclid}
\iint_{A\times A}p_t(x,y)\mathcal{H}^{2-2q}(dx)\mathcal{H}^{2-2q}(dy)\asymp t^{-q}\quad \text{ as }t\to 0.
\end{equation}
where $a \asymp b$ means that the ratio $ \ln a/ \ln b$ converges to $1$. This elementary Euclidean argument can be formulated in the LQG context. Consider the Liouville heat kernel $\p^\gamma_t(x,y)$ of the Brownian motion associated to the metric tensor \eqref{tensor}, which has been constructed in \cite{GRV2} (see also \cite{MRVZ} for further properties). A moment of thought shows that finding the quantum exponent $\Delta$ of the set $A$ consists in finding the number $ \Delta$
 such that
\begin{equation}
\E^X[\iint_{A\times A}\p^\gamma_t(x,y)\,e^{\gamma(1- \Delta)(X(x) + X(y))}\mathcal{H}^{2-2q}(dx)\mathcal{H}^{2-2q}(dy)]\asymp t^{- \Delta}\quad \text{ as }t\to 0.
\end{equation} 
The authors in  \cite{David-KPZ} claim that $\Delta$ is related to $q$ by the KPZ formula \eqref{kpz_scaling}. Their argument is based on the analysis of the metric through the {\bf Mellin--Barnes transform}. 
To explain how this comes about, it is believed that the   Liouville heat kernel should have the following behaviour:   
\begin{equation}\label{representationheat}
 \frac{c}{t} \exp \big(- \frac{{\bf d}_\gamma(x,y)^{\frac {\beta}{\beta-1}}}{ c t^{\frac{1}{\beta-1}}}\big) \leq \p^\gamma_t(x,y) \leq  \frac{C}{t} \exp
\big(-  \frac{{\bf d}_\gamma(x,y)^{\frac {\beta}{\beta-1}}}{C t^{\frac{1}{\beta-1}}}\big)\quad \text{ as } t \to 0 ,\end{equation} 
where $C,c>0$ are two positive constants, ${\bf d}_\gamma(x,y)$ is the  distance associated to the (formal) tensor \eqref{tensor} and where the exponent $\beta=\beta(\gamma)$ should represent the Hausdorff dimension associated to the (conjectural) metric space $(\R^2,{\bf d}_\gamma)$.  
The idea   is then to claim that if relation \eqref{representationheat} holds, then one can extract the metric ${\bf d}_\gamma$ from the heat kernel by using the Mellin--Barnes transform  
\begin{equation*}  
\M^\gamma_s(x,y) = \int_0^{\infty} t^{-s}\p^\gamma_t(x,y) dt.
 \end{equation*}
Indeed, assuming ~\eqref{representationheat}, a straightforward computation (using the change of variable $t=u \, {\bf d}_\gamma(x,y)^{\beta}$)  gives that for any $s \in ]0,1[$, the Mellin--Barnes transform should behave as (when ${\bf d}_\gamma(x,y) \to 0$) 
\begin{equation}\label{heuristicDavid}
\M^\gamma_s(x,y) \asymp    \frac{1}{{\bf d}_\gamma(x,y)^{\beta \, s}} \,.\end{equation}
In particular we emphasise that for $\gamma = 0$, this is an exact relation: we get
$$
M^0_s(x,y) = C \frac1{|x-y|^{2s}} 
$$
where $C = \int_0^\infty u^{s-1} e^{-1/(2u)} du< \infty$. This yields an exact connection between the Mellin--Barnes transform and the Euclidean distance and so allows one to compute the Hausdorff dimension from the Mellin-Barnes transform.
 
Of course, such a discussion is heuristic, as even the distance ${\bf d}_\gamma(x,y)$ remains out of reach so far and it is not clear that the heat kernel has the shape \eqref{representationheat} (see \cite{MRVZ} for a thorough discussion as well as theorems in this direction), except of course in the Euclidean case when $\gamma=0$.  However, we point out that the analogy with the Euclidean case $\gamma = 0$ has limitations and hides some more complex issues.
For instance it is not entirely clear what is the regime in which \eqref{representationheat} should be valid. In fact, it is \emph{a priori} natural to suppose that \eqref{representationheat} should hold when $x$ and $y$ are two fixed points and $t$ goes to $0$. On the other hand, the heat kernel KPZ formula that we will prove in this paper reflects the fact that the Mellin--Barnes transform is mainly dominated by another  regime of values of $t$, i.e. $t$ goes to $0$ simultaneously with $|y-x|$. 
Actually, we will show that the main contribution to the Mellin--Barnes transform comes from diffusive scales, i.e. those such that $t\approx |y-x|^2$ (a main technical difficulty of this article is to show that the Mellin--Barnes transform is not strongly affected by the short time scales $t\ll  |y-x|^2$, when the associated Brownian bridge in formula \eqref{dec:MBT} below concentrates along length minimizing paths and relates in a more complex fashion to the ambient Liouville volume measure).

The main goal of this paper is to prove rigorously this heat kernel derivation of the KPZ formula, thereby yielding a geometrically more intrinsic formula for Liouville quantum gravity. 
Our approach is as follows. Given a measurable set $A \subset \R^2$, and $s\in (0,1)$, we wish to define its (quantum) $s$-capacity as
\begin{align*}
C^{\gamma}_s(A)= \sup\Big\{\Big(\int_{A\times A} \M^\gamma_s(x,y)\,\mu(dx)\mu(dy)\Big)^{-1} \Big\}.
\end{align*}
where the supremum runs over all the Borel probability measures $\mu$ supported by $A$. (In fact, for technical reasons we will work with the heat kernel of the Liouville Brownian motion killed at rate $\alpha>0$, though this change is unimportant; see \eqref{MB1} and \eqref{CapMB}).  Then we define the quantum dimension of $A$ as
\begin{equation}\label{dimcap}
{\rm dim}_{\gamma}(A)=\sup\{s> 0;C^\gamma_s(A)>0\}=\inf\{s\geq 0;C^\gamma_s(A)=0\}.
\end{equation}
Note that in making this definition, we do not take advantage of the Euclidean structure which exists on $\R^2$, in as much as the Liouville measure, the Liouville Brownian motion, and therefore its heat kernel, are universal objects. One should however point out that in restricting the supremum defining the quantum capacity of a set to Borel measures, we are assuming that the Liouville metric has the same topology as the ambient Euclidean space. This might seem like a nontrivial assumption, however it is one that is commonly admitted and is in fact proved in a sense for the scaling limit of random $p$-angulations \cite{LG}. 
  
Note also that if indeed \eqref{representationheat} holds, then \eqref{heuristicDavid} shows that \eqref{dimcap} is just a multiple of the standard Hausdorff dimension in the metric space $(\R^2, {\bf d}_\gamma)$.
In particular when $\gamma=0$, it can be seen that ${\rm dim}_{0}(A)$ defined by \eqref{dimcap} is just one-half of the standard Hausdorff dimension of $A$ in terms of the Euclidean distance. Our main result (whose precise version can be found in Theorem \ref{T}) is then the following 
\begin{theorem}\label{T:intro}
Let $A\subset \R^2$ be a fixed bounded subset,  then almost surely in the the free field $X$, we have the KPZ formula
$${\rm dim}_{0}(A)= \big(1+\frac{\gamma^2}{4}\big){\rm dim}_{\gamma}(A)-\frac{\gamma^2}{4}{\rm dim}_{\gamma}(A)^2.$$ 
 \end{theorem}

 \begin{remark}
The KPZ formula is sometimes stated in terms of the Euclidean and quantum scaling exponents (resp. $q$ and $\Delta$). These are given by $q = 1- {\rm dim}_{0}(A)$ and $\Delta = 1- {\rm dim}_\gamma(A)$ respectively, and the above relation becomes
\begin{equation}\label{kpz_scaling}
q = (1-\frac{\gamma^2}4)\Delta + \frac{\gamma^2}{4} \Delta^2.
\end{equation}
\end{remark}

\begin{remark} 
It is important to note here that the subset $A\subset \R^2$ should not depend on the free field $X$, otherwise the KPZ formula does not hold anymore. This can be seen by considering the set of $\gamma$-thick points (further examples are investigated  in \cite{Aru} for the Hausdorff and Minkowski versions of KPZ).
\end{remark}

\section{Setup}\label{sec.setup}

We now discuss in more details the setup and the assumptions for Theorem \ref{T:intro}. We will work with a whole plane massive free field but our techniques may be adapted to other setup as well (massless Gaussian free field on bounded domains for instance).

\subsection{Massive Free Field}    
We consider a whole plane Massive Gaussian Free Field (MFF) $X$ with mass $m$ (see \cite{glimm,She07} for an overview of the construction of the MFF).   This is a centered Gaussian distribution in the sense of Schwartz with covariance kernel  given by 
\begin{equation}\label{MFF1}
\forall x,y \in \R^2,\quad G_m(x,y)=\int_0^{\infty}e^{-\frac{m^2}{2}u-\frac{|x-y|^2}{2u}}\frac{du}{2 u}=\ln_+\frac{1}{|x-y|}+g_m(x,y).
\end{equation}
for some continuous and bounded function $g_m$, which decays exponentially fast to $0$ when $|x-y|\to\infty$ (recall that $\ln_+(x)= \max(0,\ln x)$ for $x>0$). We denote by $\Pb^X$ and $\E^X$ the law and expectation with respect to this free field.

\subsection{Liouville measure and Liouville Brownian motion}\label{LMLBM}
We consider a coupling constant $\gamma \in [0,2[$ and consider the formal metric tensor 
$$g=e^{\gamma X(x)-\frac{\gamma^2}{2}\E[X^2]}\,dx^2.$$
The volume form of this metric tensor is  a Gaussian multiplicative chaos \cite{cf:Kah,review} with respect to the Lebesgue measure $dx$ 
\begin{equation}
M_\gamma(dx)=e^{ \gamma X(x)-\frac{\gamma^2}{2}\E[X(x)^2]}\,dx.
\end{equation}
We also consider the associated  Liouville Brownian Motion (LBM for short). More precisely, we consider a planar Brownian motion $B$. We assume that the Brownian motion and the free field $X$ are constructed on the same probability space. We denote by $\Pmp{\Bvec}{x}$ and $\Emp{\Bvec}{x}$ the probability law and expectation of this Brownian motion when starting from $x$. We will also consider the annealed probability laws $\Pb_x=\Pb^X\otimes \Pmp{\Bvec}{x}$ and the corresponding expectation $\E_x$. 
$\Pb^X$-almost surely, we consider the unique Positive Continuous Additive Functional (PCAF) $F$ associated to the measure $M_\gamma$, which is defined under $\Pmp{\Bvec}{x}$ for all starting point $x\in\R^2$ (see  \cite[section 2.4]{GRV1}).
Then, $\Pb^X$-almost surely, the law of the LBM  under $\Pmp{\Bvec}{x}$  is given by
$$ \LB_t=\Bvec_{F(t)^{-1}} $$
for all $x\in\R^2$. Furthermore this PCAF can be understood as a Gaussian multiplicative chaos with respect to the occupation measure of the Brownian motion $\Bvec$ 
\begin{equation}\label{def:F}
F(t)=\int_0^t e^{ \gamma X(\Bvec_r)-\frac{\gamma^2}{2}\E^X[X^2(\Bvec_r)]}\,dr.
\end{equation}
The LBM is a  Feller Markov process with continuous sample paths \cite[section 2.7]{GRV1}.    Let us denote by $(P^\gamma_t)_{t\geq 0}$ the associated Liouville semigroup (which is random as it depends on $X$). $\Pb^X$-almost surely, this semigroup is absolutely continuous  with respect to the Liouville measure $M_\gamma$  and there exists a measurable function $ \p^\gamma_t(x,y)$, called Liouville heat kernel such that for all $x\in\R^2$ and any measurable bounded function $f$ (see \cite[section 2]{GRV2}) 
\begin{equation}
P^\gamma_tf(x)=\int_{\R^2}f(y)\p^\gamma_t(x,y)\,M_\gamma(dy).
\end{equation}
Actually, one can even show that $\p^\gamma_t(x,y)$ is a continuous function of $(t,x,y)$  (see \cite[section 3]{MRVZ}).
In what follows, we will also consider the standard heat kernel $p_t(x,y)$  of the planar Brownian motion $B$ on $\R^2$.

\subsection{ Mellin-Barnes transform}
The Mellin-Barnes transform of the Liouville heat kernel is defined for $s\geq 0$ by
\begin{equation}\label{MB11}
\int_{0}^\infty  t^{-s}    \p^\gamma_t(x,y) dt .
\end{equation}
Actually, we will not work with this definition because it involves the large $t$ behaviour of the heat kernel, which does not play a significant part for the definition of quantum capacity dimension below or the KPZ relation but raises additional technical difficulties when working on the whole plane. So we will reduce the effect of the large $t$ behaviour of the heat kernel by considering instead the following Mellin-Barnes like transform for $\alpha>0$
\begin{equation}\label{MB1}
\M^{\alpha, \gamma}_s( x,y) = \int_{0}^\infty  t^{-s} e^{-\alpha t}    \p^\gamma_t(x,y) dt,
\end{equation}
which we call $(\alpha,\gamma)$-type Mellin-Barnes  transform, or $(\alpha,\gamma)$-MBT for short.


\subsection{Capacity dimension}    

 Given a Borel set $A\subset \R^2$, we define the quantum $s$-capacity of $A$ by
\begin{align}
C^{\alpha,\gamma}_s(A)=&\sup\Big\{\Big(\int_{A\times A}\M^{\alpha,\gamma}_s(x,y)\mu(dx)\mu(dy)\Big)^{-1} \Big\} \label{CapMB}
 \end{align}
where the supremum runs over the Borel  probability measures $\mu$ such that $\mu(A)=1$.
The quantum capacity dimension of the  set $A$ is defined as the quantity
\begin{equation}
{\rm dim}_{\gamma}(A)=\sup\{s\geq 0;C^{\alpha,\gamma}_s(A)>0\}=\inf\{s\geq 0;C^{\alpha,\gamma}_s(A)=0\}. \label{MBcap2}
\end{equation}
A priori, this definition depends on the killing rate $\alpha$. However, we will see that, at least when $A$ is a bounded set, this quantity does not depend on $\alpha$. 

When $\gamma = 0$, note that 
\begin{align*}
\M_s^{\alpha, 0 } (x,y) &= \int_0^\infty t^{-s} e^{ - \frac{|x-y|^2}{2t}} \frac1{2\pi t } e^{- \alpha t } dt \\
& =  C |x-y|^{-2s} \int_0^\infty u^{-s-1} e^{- \frac1{2u}} e^{-\alpha |x-y|^2 u } du\\
&\le C |x-y|^{-2s}
\end{align*}
for some constant $C$. A similar lower bound can be established, with a different constant $C$, provided that $A$ is bounded, say $A \subset B(0,1/2)$, which we assume without loss of generality henceforth. Indeed, we get
\begin{align*}
M_s^{\alpha, 0 } (x,y )& \ge C |x-y|^{-2s} \int_0^{\frac1{|x-y|^2}} u^{-s-1} e^{- \frac1{2u}} du \\
& \ge  C |x-y|^{-2s} \int_0^{1/4} u^{-s-1} e^{- \frac1{2u}} du = c|x-y|^{-2s}.
\end{align*}
Therefore, we deduce that 
$$
{\rm dim}_{0}(A) = \frac12 \dim_{\text{Hausdorff}, \R^2 }(A)
$$
where $\dim_{\text{Hausdorff}, \R^2 }(A)$ denotes the regular Hausdorff dimension of $A$ with respect to the standard Euclidean metric on $\R^2$.

\subsection{KPZ formula}
With these definitions we are ready to give a precise version of Theorem \ref{T:intro}.
\begin{theorem}\label{T}
Let $A\subset \R^2$ be bounded. Then ${\rm dim}_\gamma(A)$ does not depend on $\alpha$. Furthermore, the Euclidean and quantum capacity dimension of $A$ are $\Pb^X$-almost surely related to each other by the KPZ formula
\begin{equation}
{\rm dim}_{0}(A)= \big(1+\frac{\gamma^2}{4}\big){\rm dim}_{\gamma}(A)-\frac{\gamma^2}{4}{\rm dim}_{\gamma}(A)^2.\label{kpz}
\end{equation}
\end{theorem}

\section{Proof of the heat kernel KPZ formula}    
In this section, $C$ will stand for a constant that does not depend on any relevant quantity. It may change along lines in the computations. We define 
\begin{equation}\label{xi}
\xi(q)=(2+\frac{\gamma^2}{2})q-\frac{\gamma^2}{2}q^2
\end{equation}
 for $q>0$:  this function coincides with the power law spectrum of the measure $M_\gamma$ (see \cite{review} for more on this). Let us denote by $\mathcal{P}$ the set of Borel probability measures  supported by $A$. Without loss of generality, we may assume that $A\subset B(0,1/2)$. 

\subsection{Background on the law of the Brownian bridge} 
We denote by  $\Pbr{\Bvec}{x}{y}{t}$ and $\Ebr{\Bvec}{x}{y}{t}$ the law and expectation of the Brownian bridge from $x$ to $y$ with lifetime $t$. Therefore, under $\Pbr{\Bvec}{x}{y}{t}$, $(\Bvec_s)_{0\leq s \leq t}$ is a Brownian bridge from $x$ to $y$ with lifetime $t$.

Recall the basic lemma (see (6.28) in \cite{KaratzasShreve})
\begin{lemma}\label{absBB}
We have the following absolute continuity relation for the Brownian bridge for $s<t$:
\begin{equation*}
\Ebr{\Bvec}{x}{y}{t}[  G( (\Bvec_u)_{u \leq s}   )   ]  = \Emp{\Bvec}{x}[   G( (\Bvec_u)_{u \leq s}   )  \frac{t}{t-s}e^{\frac{|y-x|^2}{2 t}-\frac{|\Bvec_s-x|^2}{2 (t-s)}  }     ].
\end{equation*}
In particular, if $(\Bvec_s, 0 \le s \le t)$ has law $\Ebr{\Bvec}{x}{y}{t}$, then  $(\lambda^{-1} \Bvec_{\lambda^2 s}, 0 \le s \le t /\lambda^2)$ has law $\Ebr{\Bvec}{x/\lambda}{y/\lambda}{t/\lambda^2}$.
\end{lemma}

We will also need the following elementary representation of a Brownian bridge (see, e.g., (6.29) in \cite{KaratzasShreve}). 
\begin{lemma}\label{BBlinear}
Let $\Bvec$ be a standard Brownian motion and let $t\ge0$, and $x, y \in \R^2$. Then $$b_s = x + \Bvec_s - \frac{s}{t}  (y - x - \Bvec_t), \ \  0 \le s \le t$$ defines a Brownian bridge from $x$ to $y$ of duration $t$.
\end{lemma}

\subsection{Brownian Bridge decomposition}

We recall the following result
\begin{theorem}[see \cite{spectral}]\label{generalformula}
$\Pb^X$ almost surely, for each $x,y\in\R^2$ and any continuous function $G:\R_+\to\R_+$
\begin{equation} \label{transformBB}
\int_{0}^\infty  G(t)    \p^\gamma_t(x,y) dt =   \int_0^\infty  \Ebr{\Bvec}{x}{y}{t}[ G(  F(t))  ]  p_t(x,y) dt  .
\end{equation}
\end{theorem}

By applying Theorem \ref{generalformula} with the function $G(t)=t^{-s} e^{-\alpha t} $ for $\alpha,s>0$ we get the following   Brownian bridge decomposition of the Mellin--Barnes transform. 
\begin{equation}
\label{dec:MBT}
\M^{\alpha, \gamma}_s( x,y) = \int_{0}^\infty  t^{-s} e^{-\alpha t}    \p^\gamma_t(x,y) dt=\int_0^\infty \Ebr{\Bvec}{x}{y}{t}[e^{-\alpha F(t)}F(t)^{-s} ]\frac{e^{-\frac{|y-x|^2}{2t}}}{2\pi t}\,dt.
\end{equation}
Likewise,
$$
C_s^{\alpha, \gamma}(A) =  \sup\Big\{\Big(\int_{A\times A}\int_0^\infty \Ebr{\Bvec}{x}{y}{t}[e^{-\alpha F(t)}F(t)^{-s} ] p_t(x,y)\,dt\mu(dx)\mu(dy)\Big)^{-1}\Big\},
$$

We will use this relation throughout the paper.

\subsection{Background on the measure based KPZ formula} \label{backKPZ}

Here we recall a few basic facts about the measure based KPZ formula as stated in \cite{Rnew10,Rnew4}. Yet, we will need to extend slightly the framework. Given an atom free Radon measure $\mu$ on $\R^2$ and $s\in [0,1]$, we define for $\delta>0$ and $a\geq 1$
\begin{equation}\label{hauss}
H^{s,\delta,a}_\mu(A)= \inf \big\{\sum_k \mu(B(x_k,a r_k))^{s} \big\}
\end{equation}
where the infimum runs over all the covering $(B(x_k,r_k))_k$ of $A$ with closed Euclidean balls with radius $r_k\leq \delta$. Notice the factor $a$, which differs from \cite{Rnew10,Rnew4}.  
The mapping $\delta>0\mapsto H^{s,\delta,a}_\mu(A)$ is decreasing and we can define:
$$H^{s,a}_{\mu}(A)=\lim_{\delta\to 0}H^{s,\delta,a}_{\mu}(A).$$
$H^{s,a}_\mu$ is a metric outer measure on $\R^2$ its restriction to the $\sigma$-field of $H^{s,a}_\mu$-measurable sets, which contains all the  Borel sets, is a measure. The $\mu$-Hausdorff dimension of the set $A$ is then defined as the value
\begin{equation} 
{\rm dim}_{\mu,a}(A)=\inf\{s\geq 0; \,\,H^{s,a}_\mu(A)=0\}=\sup\{s\geq 0; \,\,H^{s,a}_\mu(A)=+\infty\}.
\end{equation}
Notice that ${\rm dim}_{\mu,a}(A)\in [0,1]$.
 
\begin{remark}
When $\mu$ is the Lebesgue measure, the notion of Hausdorff dimension does not depend on $a$ so that we write ${\rm dim}_{Leb}(A)$ instead of ${\rm dim}_{Leb,a}(A)$. This comes from the fact that ${\rm  Leb}(r A)=r^2{\rm Leb}(A)$ for all Borel set $A$ and $r>0$. The measure $M_\gamma$ is far from having such nice scaling properties. In particular, the readers that are used to "smooth" measures (meaning satisfying the doubling condition) may find part of our forthcoming proofs unnatural. So we stress that the Liouville measure  is far from satisfying the doubling condition.  We can show that we have (see appendix \ref{sec:doubling})
$$\underset{r \to 0}{\overline{\lim}}\sup_{x\in B(0,1)}\frac{M_\gamma(B(x,2r))}{M_\gamma(B(x,r))^{1-\eta}} 
<+\infty, 
\text{ if }\eta>\frac{\gamma^2}{4+\gamma^2} 
$$
and we believe that if $\eta<\frac{\gamma^2}{4+\gamma^2}$ then the corresponding limit is infinite. A proof of the analogue statement is in fact provided in the appendix in the case of one-dimensional lognormal multiplicative cascades.
\end{remark}

We fix $a\geq 1$. In what follows, given a compact set $A$ of $\R^2$, we define its Hausdorff dimensions 
${\rm dim}_{Leb}(A)$, ${\rm dim}_{M_\gamma,a}(A)$   computed as indicated above with $\mu$ respectively equal to the Lebesgue measure or $M_\gamma$.

\begin{theorem}\label{KPZmeasure}{{\bf Measure-based KPZ formula.}}
Fix $a\geq 1$. Let $A$ be a compact set of $\R^2$. $\Pb^X$ almost surely, we have the relation 
$${\rm dim}_{Leb}(A)=(1+\frac{\gamma^2}{4}){\rm dim}_{M_\gamma,a}(A)  - \frac{\gamma^2}{4}{\rm dim}_{M_\gamma,a}(A)^2.$$ 
In particular, the quantity ${\rm dim}_{M_\gamma,a}(A)$ does not depend on $a\geq 1$ (and thus we skip the index $a$).
\end{theorem}

\begin{remark} 
It is easy to check that ${\rm dim}_{Leb}(A)={\rm dim}_{0}(A)$. This is because in the definition of Hausdorff dimension of a set $A$, it suffices to consider coverings of $A$ using Euclidean balls of diameter bounded by $\delta$ (see e.g. the discussion after (2.16) in \cite{Falconer}). 
\end{remark}

The proof of the KPZ formula in the case $a=1$ is based in \cite{Rnew10,Rnew4} on a measure based version of the Frostman lemma involving quarters of balls. 
The case $a>1$ in the  definition \eqref{hauss} allows us more flexibility in the shapes of the sets used in this modified Frostman lemma. Thus we are going to state this lemma.

\begin{remark}
All the results below work also in the case when, instead of the ball $B(\frac{x+y}{2},\frac{(a-1)}{2}|y-x|)$, we consider the ``tube" $C(x,y,a)$  made up of those points located at (Euclidean) distance less than $(a-1)\frac{|x-y|}{2}$ of the segment $[x,y]$, namely 
$$C(x,y,a)=\{z\in\R^2; {\rm dist}(z,[x,y])\leq (a-1)\frac{|x-y|}{2}\}.$$
\end{remark}
 
\begin{lemma}[Modified Frostman lemma]\label{lem:frost}
Assume $a>1$ and that $\mu,\nu$ are two Radon measures on $\R^2$. Assume further that $\nu$ is a probability measure supported by a compact set $A\subset B(0,1)$. If for $q\in [0,1]$
\begin{equation}\label{eq:frost}
\int_A\int_A\frac{\nu(dx)\nu(dy)}{ \mu\Big(B(\frac{x+y}{2},\frac{(a-1)}{2}|y-x|)\Big)^q}<+\infty
\end{equation}
then ${\rm dim}_{\mu,a}(A)\geq q$.
\end{lemma}

The proof of the measured-based KPZ formula (Theorem \ref{KPZmeasure}) can then be completed exactly as in \cite{Rnew4,Rnew10} with the help of Lemma \ref{lem:frost}.

\medskip 
\noindent {\it Proof of Lemma \ref{lem:frost}.} Let us define the function:
$$\forall x\in A,\quad g(x)=\int_K\frac{\nu(dy)}{\mu\Big(B(\frac{x+y}{2},\frac{(a-1)}{2}|y-x|)\Big)^q}.$$
Observe that the assumptions imply that $\int_Ag(x)\nu(dx)< \infty$. We deduce that 
$$\nu\big(\{x\in A;g(x)\leq L\}\big)\to 1,\quad \text{ as }L\to\infty.$$ Therefore we can find $L$ large enough such that the set $A_L=\{x\in A;g(x)\leq L\}$ satisfies $\nu(A_L)\geq \frac{1}{2}$. Let us consider a covering $(B(x_n,r_n))_n$ of $A$ with balls of radius less than $\delta$. We consider the subsequence $(B(x_{n_k},r_{n_k}))_{n_k}$ of balls that intersect $A_L$. It is obvious that this subsequence forms a covering of $A_L$. For each $n_k$, there exists $y_{n_k}$ in $A_L \cap B(x_{n_k},r_{n_k})$. For all $y \in B(x_{n_k},r_{n_k})$, the set $B(\frac{y_{n_k}+y}{2},\frac{(a-1)}{2}|y-y_{n_k}|) $ is contained in $B(x_{n_k},ar_{n_k})$. Hence
\begin{equation*}
 \mu\Big(B(\frac{y_{n_k}+y}{2},\frac{(a-1)}{2}|y-y_{n_k}|) \Big)^q  \leq \mu( B(x_{n_k},ar_{n_k})  )^q.
\end{equation*}
Therefore, we get, passing to quotients and integrating over $B(x_{n_k}, r_{n_k})$ with respect to $\nu$ (and noting that the integral over the ball is dominated by the integral over $A$ since $\nu$ is supported on $A$),
\begin{equation*}
\frac{\nu(  B(x_{n_k},r_{n_k}) )}{\mu(  B(x_{n_k},ar_{n_k}) )^q  } \leq  \int_{A}\frac{\nu(dy)}{ \mu\Big( B(\frac{y_{n_k}+y}{2},\frac{(a-1)}{2}|y-y_{n_k}|) \Big)^q}   = g(y_{n_k})  \leq L.
\end{equation*} 
This leads to 
\begin{equation*}
\sum_k  \mu(  B(x_{n_k},ar_{n_k}) )^q  \geq \frac{1}{L} \sum_k  \nu(  B(x_{n_k},r_{n_k}) ) \geq \frac{1}{L} \nu(  A_L)  \geq \frac{1}{2L}.
\end{equation*}
Thus $\sum_n  \mu(  B(x_{n},ar_{n}) )^q  \ge 1/(2L)$. Since this covering of $A$ was arbitrary, we get the desired result. \qed

\subsection{Upper bound} 
\begin{rem}
For the time being, we have not proved yet that the notion of quantum capacity dimension ${\rm dim}_{\gamma}(A)$ does not depend on the exponent $\alpha$ in \eqref{MBcap2}. Therefore, we use the notation ${\rm dim}_{\alpha,\gamma}(A)$ instead of ${\rm dim}_{\gamma}(A)$ to keep track of this "a priori" dependence until we have proved this statement.
\end{rem} 

We first establish the bound 
\begin{equation}\label{UB}
{\rm dim}_{cap}(A)\geq \big(1+\frac{\gamma^2}{4}\big){\rm dim}_{\alpha,\gamma}(A)-\frac{\gamma^2}{4}{\rm dim}_{\alpha,\gamma}(A)^2.
\end{equation}
The proof of this bound is based on the following two lemmas:
\begin{lemma}\label{minoration}
$\Pb^X$ a.s., there exists a constant $D$ (random and measurable with respect to $X$) such that for all $x,y\in B(0,1/2)$, $s\in[0,1]$ 
\begin{equation}\label{E:UB}
\int_0^\infty  \Ebr{\Bvec}{x}{y}{t} [  e^{-\alpha F(t)}F(t)^{-s} ] p_t(x,y)\,dt\geq  \frac{D}{(\mu(x,|y-x|) +\mu(y,|y-x|))^s},
\end{equation}
where $(x,r)\mapsto\mu(x,r)$ is defined for $x\in \R^2$ and $r>0$ by
$$\mu(x,r):=\int_{B(x,2r)}(1+\ln\frac{1}{|x-z|})\,M_\gamma(dz).$$  
\end{lemma}

\begin{lemma}\label{capM}
$\Pb^X$ a.s., for all $s\geq 0$, if one can find a probability measure $\nu$ on $A$ such that 
\begin{equation}\label{lem:sup}
 \int_{A\times A}  \frac{1}{(\mu(x,|y-x|)+\mu(y,|y-x|))^s}\,\nu(dx)\nu(dy) <\infty
\end{equation} 
 then $\xi(s)/2\leq {\rm dim}_{cap}(A)$.
\end{lemma}

\medskip
Let us admit for a while the two above lemmas and complete the proof of \eqref{UB}. The strategy will be to show (using Lemma \ref{muM}) that \eqref{lem:sup} basically entails 
$$ \int_{A\times A}  \frac{1}{M_\gamma(B(\frac{x+y} 2,4|y-x|))^s}\,\nu(dx)\nu(dy)< \infty$$
 and then use Lemma \ref{lem:frost} to get the lower bound for the capacity dimension with $a=9$. 

In the case where ${\rm dim}_{\alpha, \gamma}(A) = 0$ there is nothing to verify so we may assume without loss of generality that ${\rm dim}_{\alpha, \gamma}(A) >0$. Thus let us consider $s\geq 0$ such that  $s<{\rm dim}_{\alpha,\gamma}(A)$, that is $C_s^{\alpha,\gamma}(A)>0$. One can thus find a probability measure $\nu$ on $A$ which is such that 

\begin{equation*}
\int_{A\times A} \Big( \int_0^\infty  \Ebr{\Bvec}{x}{y}{t} [  e^{-\alpha F(t)}F(t)^{-s} ] p_t(x,y)\,dt \Big) \nu(dx) \nu(dy)  < \infty\,.
\end{equation*}

From Lemma \ref{minoration}, we obtain that 
\begin{align*}
\int_{A\times A}   \frac{1}{(\mu(x,|y-x|)+\mu(y,|y-x|))^s} \nu(dx) \nu(dy)  < \infty\,,
\end{align*}
which in turns implies by Lemma \ref{capM} that $\xi(s)/2\leq {\rm dim}_{cap}(A)$. Since this inequality  holds for all $s< {\rm dim}_{\alpha,\gamma}(A)$ and since $q\mapsto \xi(q)$ is increasing, we obtain our desired upper bound 
$$\hspace{5cm }\frac{\xi({\rm dim}_{\alpha,\gamma}(A))}{2} \leq {\rm dim}_{cap}(A).\hspace{ 6cm }\qed$$

\medskip

\noindent {\it Proof of Lemma \ref{minoration}.} We first bound from below the $(\alpha,\gamma)$-type Mellin-Barnes transform by selecting the more relevant (i.e. diffusive) scales
\begin{align*}
\int_0^\infty  \Ebr{\Bvec}{x}{y}{t}&[  e^{-\alpha F(t)}F(t)^{-s} ] p_t(x,y)\,dt \\
\geq & \int_{|y-x|^2/2}^{|y-x|^2} \Ebr{\Bvec}{x}{y}{t} [  e^{-\alpha F(t)}F(t)^{-s} ] p_t(x,y)\,dt\\
\geq & \int_{|y-x|^2/2}^{|y-x|^2} \,\frac{dt}{2\pi e  t}\times \inf_{|y-x|^2/2\leq t\leq |y-x|^2} \Ebr{\Bvec}{x}{y}{t} [  e^{-\alpha F(t)}F(t)^{-s} ] \\
= & \frac{\ln 2}{2\pi e } \times \inf_{|y-x|^2/2\leq t\leq |y-x|^2} \Ebr{\Bvec}{x}{y}{t} [  e^{-\alpha F(t)}F(t)^{-s} ] .
\end{align*}

\begin{rem}
In what follows, the strategy is more or less the following. By using the Jensen inequality, we should have 
$$ \inf_{|y-x|^2/2\leq t\leq |y-x|^2} \Ebr{\Bvec}{x}{y}{t} [  e^{-\alpha F(t)}F(t)^{-s} ] \geq  \exp\left({-\alpha \Ebr{\Bvec}{x}{y}{t=|y-x|^2}[F(t)]} \right) \left(\Ebr{\Bvec}{x}{y}{t=|y-x|^2}[F(t)] \right)^{-s}.$$
Then we just have to compute $\Ebr{\Bvec}{x}{y}{t=|y-x|^2}[F(t)]$, which approximatively takes on the form
$$\Ebr{\Bvec}{x}{y}{t=|y-x|^2}[F(t)] \asymp \int_{\R^2} g(x,y,z)M_\gamma(dz),$$ for some function $g$ that gives almost all its mass to the ball $B(\frac{x+y}{2},|y-x|)$ (and which is related to the Green function of Brownian motion in that ball). Yet, this function $g$ does not vanish outside this ball and possesses a non trivial behaviour for large $z$. Basically, the following proof follows this idea except that we will try to get rid of this long range dependence while getting in the end a   tractable function $g$.
\end{rem}
Let us introduce the event 
$$A=\{\forall s\in [0,t]; \Bvec_s\in B(\frac{x+y}{2};|y-x|)\}$$ and
the stopping time
$$T(x,r)=\inf\{s\geq 0; \Bvec_s\not \in B(x,r)\}.$$

Now we estimate the quantity  $\Ebr{\Bvec}{x}{y}{t} [  e^{-\alpha F(t)}F(t)^{-s} ] $ for $|y-x|^2/2\leq t\leq |y-x|^2$. First we write  
\begin{align}
  \Ebr{\Bvec}{x}{y}{t} [  e^{-\alpha F(t)}F(t)^{-s} ] &\geq \Ebr{\Bvec}{x}{y}{t} [  e^{-\alpha F(t)}F(t)^{-s} \ind_{A}]\nonumber\\
  &=  \Ebr{\Bvec}{x}{y}{t} [  e^{-\alpha F(t)}F(t)^{-s} |A] \Pbr{\Bvec}{x}{y}{t} (A).\label{porc1}
  \end{align}
The important point in the following is to notice that there exists a constant $c>0$ such that for any $x,y$ and   $|y-x|^2/2\leq t\leq |y-x|^2$ 
\begin{equation}\label{porc2}
c^{-1}\leq \Pbr{\Bvec}{x}{y}{t} (A)\leq c.
\end{equation}
Indeed by translation invariance we may assume that $x = 0$. Then by applying the scaling in Lemma \ref{absBB}, \eqref{porc2} reduces to proving that 
$$
c^{-1} \leq \Pbr{\Bvec}{0}{1}{t}(A') \leq c
$$ 
where $t\in [1/2, 1]$ is arbitrary and $A'$ is the event $A' = \{\forall s\in [0,t]; \Bvec_s\in B((1/2,0);1)\}$. For a fixed $t>0$ observe that $\Pbr{\Bvec}{0}{1}{t}(A') $ is bounded away from 0 and 1, hence the desired statement follows from the continuity in $t$ of this probability, which can be seen from Lemma \ref{BBlinear}. %

Observe that the mapping $u\mapsto e^{-\alpha u}u^{-s} $ is convex. 
By applying the Jensen inequality, we deduce
\begin{align}
 \Ebr{\Bvec}{x}{y}{t} [  e^{-\alpha F(t)}F(t)^{-s} |A] \geq e^{-\alpha  \Ebr{\Bvec}{x}{y}{t} [  F(t)  |A]}( \Ebr{\Bvec}{x}{y}{t} [  F(t)  |A])^{-s}.\label{porc3}
  \end{align}
Let us compute   $\Ebr{\Bvec}{x}{y}{t} [  F(t)  |A]$. We have from \eqref{porc2} 
\begin{equation}\label{porc4}
\Ebr{\Bvec}{x}{y}{t} [  F(t)  |A]\leq c\Ebr{\Bvec}{x}{y}{t} [  F(t)  \ind_A].
\end{equation}
Then, by symmetry of the time reversed law  of the Brownian bridge,   the fact that the ball $B(\frac{x+y}{2};|y-x|)$ is contained in the ball $B(x,2|y-x|)$ and Lemma \ref{absBB}  
\begin{align}
\Ebr{\Bvec}{x}{y}{t} [  F(t)  \ind_A]\leq & \Emp{\Bvec}{x} [  F(t/2)  \ind_{\{T_{B(x,2|y-x|)}>t/2\}}2e^{\frac{|y-x|^2}{2 t}-\frac{|\Bvec_{t/2}-x|^2}{t}  }   ]\nonumber \\
& \ \ \ \  + \Emp{\Bvec}{y} [  F(t/2)  \ind_{\{T_{B(y,2|y-x|)}>t/2\}}2e^{\frac{|y-x|^2}{2 t}-\frac{|\Bvec_{t/2}-y|^2}{t}  }   ] \nonumber\\
\leq  &2 e \, \Emp{\Bvec}{x} [  F(T(x,2|y-x|))  \ind_{\{T_{B(x,2|y-x|)}>t/2\}}] + 2 e \, \Emp{\Bvec}{y} [  F(T(y,2|y-x|))  \ind_{\{T_{B(y,2|y-x|)}>t/2\}}] \nonumber\\
=& 2 e\int_{B(x,2|y-x|)}G_{B(x,2|y-x|)}(x,z)\,M_\gamma(dz) + 2 e\int_{B(y,2|y-x|)}G_{B(y,2|y-x|)}(x,z)\,M_\gamma(dz)\label{porc5}
\end{align}  
 where  $G_{B(x,2|y-x|)}(z,z')$ stands for the Green function of the Brownian motion killed upon leaving the ball 
  $B(x,2|y-x|)$. 
  It is a standard fact that we have the explicit expression
\begin{equation}\label{porc6}
G_{B(x,2|y-x|)}(x,z)=\frac{1}{\pi}\ln\frac{2|y-x|}{|x-z|}\leq \frac{1}{\pi}\ln (2{\rm diam}(A))+\frac{1}{\pi}\ln\frac{1}{|x-z|}.
\end{equation}
Let us set for $x\in\R^2$ and $r>0$
\begin{equation}\label{porc7}
\mu(x,r):= \int_{B(x,2r)}\Big( \frac{1}{\pi}\ln (2{\rm diam}(A))+\frac{1}{\pi}\ln\frac{1}{|x-z|}\Big)\,M_\gamma(dz).
\end{equation}
By gathering \eqref{porc1}+\eqref{porc2}+\eqref{porc3}+\eqref{porc4}+\eqref{porc5}+\eqref{porc6}, the fact that the mapping $u\mapsto e^{-\alpha u}u^{-s} $ is decreasing, and the fact that $A\subset B(0,1/2)$, we deduce  
\begin{align*} 
  \Ebr{\Bvec}{x}{y}{t} [  e^{-\alpha F(t)}F(t)^{-s} ]\geq & c  e^{-\alpha 4ec[\mu(x,|y-x|)+\mu(y,|y-x|)]} (4ec[\mu(x,|y-x|)+\mu(y,|y-x|)])^{-s}\\
  \geq & c  e^{-\alpha 8ec\sup_{x\in B(0,1/2)}\mu(x,1)} (4ec[\mu(x,|y-x|)+\mu(y,|y-x|)])^{-s}
\end{align*}
We complete the proof of Lemma \ref{minoration} with the following lemma to get rid of the exponential term in the above right-hand side. \qed

\begin{lemma}
There exists a random constant $C$, $\Pb^X$ almost surely finite, such that
$$\sup_{x\in B(0,1)} \mu(x,1)\leq C.$$
\end{lemma} 

\proof Recall \cite{GRV1} that for all $0<\alpha<2(1-\frac{\gamma}{2})^2$, there exists a random $C$ such that for all $x\in B(0,1)$ and $r<1$
\begin{equation}\label{modulus}
M_\gamma(B(x,r))\leq Cr^\alpha.
\end{equation}
Observe that it is enough to show that, $\Pb^X$ a.s. 
$$\sup_{x\in B(0,1)}\int_{B(x,1)}  \ln\frac{1}{|x-z|} \,M_\gamma(dz)<+\infty.$$
We have 
\begin{align*}
\int_{B(x,1)}  \ln\frac{1}{|x-z|} \,M_\gamma(dz)\leq &\sum_{n\geq 0}\int_{z\in B(x,1),2^{-n-1}|x-z|\leq 2^{-n}}  \ln\frac{1}{|x-z|} \,M_\gamma(dz)\\
\leq &\sum_{n\geq 0}  (n+1) \ln2 \,M_\gamma(B(x,2^{-n}))\\
\leq &C\sum_{n\geq 0} 2^{-\alpha n} (n+1) \ln2 .
\end{align*}
This latter quantity is obviously finite and does not depend on $x$.\qed

\medskip 

As explained above, the proof of Lemma \ref{capM} requires to compare $\mu(x,r)$ with $M_\gamma(B(x,2r))$. This is the purpose of the following Lemma.  
\begin{lemma}\label{muM}
For each $\epsilon>0$, there exists a random constant $C_\epsilon=C_\epsilon(X)<\infty$ such that $\Pb^X$ a.s., for all $x\in B(0,1)$ and $r<1$  
\begin{equation}\label{modmu}
\mu(x,r)\leq C_\epsilon M_\gamma(B(x,2r))^{1-\epsilon}.
\end{equation}
In particular, for any $x,y \in \R^2$,
\begin{align}\label{modmu2}
\mu(x,|x-y|)+\mu(y,|x-y|)& \leq 2 C_\epsilon \big(M_\gamma(B(x,2|x-y|))+M_\gamma(B(y,2|x-y|))\big)^{1-\epsilon} \nonumber \\
& \leq 4 C_\epsilon M_\gamma(B(\frac{x+y}{2},4|x-y|))^{1-\eps}\,.
\end{align}
\end{lemma}

\noindent {\it Proof of Lemma \eqref{muM}. }To see this, it is enough to prove this relation with $\bar{\mu}(x,r)= \int_{B(x,2r)}  \ln\frac{1}{|x-z|} \,M_\gamma(dz)$ instead of $\mu$. 
Also, recall  \eqref{modulus} and choose $0<\delta<2(1-\gamma/2)^2$. $\Pb^X$ a.s.,  there exists a random constant $C$ such that for all $x\in B(0,1)$ and $u<1$
\begin{equation}\label{modulus2}
M_\gamma(B(x,u))\leq Cu^\delta.
\end{equation}
Then we have
\begin{align*}
\bar{\mu}(x,r)=&\sum_{n\geq -\log_2  r }\int_{2^{-n}\leq |z-x|\leq 2^{-(n-1)}}  \ln\frac{1}{|x-z|} \,M_\gamma(dz)\\
\leq &\sum_{n\geq -\log_2 r}   n\ln2 \,M_\gamma(B(x,2^{-(n-1)}))\\
\leq & \sum_{n\geq -\log_2 r}   n\ln2 \,M_\gamma(B(x,2^{-(n-1)}))^{1-\epsilon} M_\gamma(B(x,2^{-(n-1)}))^\epsilon\\
\leq &M_\gamma(B(x,2r))^{1-\epsilon}C^\epsilon\sum_{n\geq -\log_2 r}   n\ln2 \,  2^{-\delta(n-1)\epsilon}.
\end{align*}
 The latter series converges and can be bounded independently of $x,r$ which ends the proof. \qed

\medskip
\medskip
\noindent {\it Proof of Lemma \eqref{capM}. }

\medskip
From \eqref{modmu2}, we deduce
\begin{align*}
\sup_{\nu\in\mathcal{P}} \Big( \int_{A\times A}  \frac{1}{(\mu(x,|y-x|)+\mu(y,|y-x|))^s}\,\nu(dx)\nu(dy)\Big)^{-1} & \\
& \hskip -4 cm \leq  C_\epsilon^s\sup_{\nu\in\mathcal{P}}\Big( \int_{A\times A}  \frac{1}{M_\gamma(\frac {x+y} 2,4|y-x|)^{s(1-\epsilon)}}\,\nu(dx)\nu(dy) \Big)^{-1}.  
\end{align*}
In particular, if the left-hand side is infinite, so is the right-hand side for all $\epsilon>0$. Applying Lemma \ref{lem:frost} with $a=9$, we deduce that for all $\epsilon>0$, ${\rm dim}_{M_\gamma,a}(A)\geq s(1-\epsilon)$. We deduce that ${\rm dim}_{M_\gamma,a}(A)\geq s$, then ${\rm dim}_{Leb}(A)\geq \xi(s)$ from Theorem \ref{KPZmeasure} and finally ${\rm dim}_{cap}(A)\geq \xi(s)$ because of the standard fact that ${\rm dim}_{Leb}(A)={\rm dim}_{cap}(A)$.\qed

\subsection{Lower bound.} 

Here we prove the remaining estimate
\begin{equation}\label{LB}
{\rm dim}_{cap}(A)\leq \big(1+\frac{\gamma^2}{4}\big){\rm dim}_{\alpha,\gamma}(A)-\frac{\gamma^2}{4}{\rm dim}_{\alpha,\gamma}(A)^2.
\end{equation}
Let us consider $s\geq 0$ such that $\xi(s)/2<{\rm dim}_{cap}(A)$. 
This means that we can find a Borel probability measure $\nu$ supported by $A$ and a number $\delta>0$ such that
\begin{equation}\label{eucap0}
\int_{A\times A}\frac{1}{|x-y|^{\xi(s) + 2\delta}}\nu(dx)\nu(dy)<\infty.
\end{equation}
 We define the tilted measure
$$\widetilde{\nu}(dx)=e^{s\gamma X(x)-\frac{\gamma^2s^2}{2}\E^X[X^2(x)]}\nu(dx),$$
which is non trivial (see \cite{cf:Kah,review}) and almost surely supported by $A$.

It is enough to show that the following expectation is finite
\begin{equation}\label{LQGcap}
\E^X\Big[ \int_{A\times A}\int_0^\infty t^{-s}e^{-\alpha t}\p_t^\gamma(x,y) \tilde \nu(dx) \tilde \nu(dy)\Big]<+\infty
\end{equation}
because this entails that   $\Pb^X$ a.s.  $s\leq {\rm dim}_{\alpha,\gamma}(A)$, which completes the proof of \eqref{LB}.\qed

\medskip
\noindent {\it Proof of \eqref{LQGcap}.} We have by \eqref{generalformula},
\begin{align*}
\E^X\Big[&\int_{A\times A}\int_0^\infty t^{-s}e^{-\alpha t}\p_t^\gamma(x,y)\widetilde{\nu}(dx)\widetilde{\nu}(dy)\Big]\\
=&\E^X\Big[\int_{A\times A}\int_0^\infty \Ebr{\Bvec}{x}{y}{t}[e^{-\alpha F(t)}F(t)^{-s} ] p_t(x,y)\,dt\,\widetilde{\nu}(dx)\widetilde{\nu}(dy)\Big].
\end{align*}
Now we will split the above integral over time in two parts to analyze its contribution depending on the scales  $|x-y|^2\leq t$ or $|x-y|^2\geq t$.
 
\subsubsection*{Short scales}
 In this subsection, we compute
 $$\E^X\Big[\int_{A\times A}\int_0^{|x-y|^2} \Ebr{\Bvec}{x}{y}{t}[e^{-\alpha F(t)}F(t)^{-s} ] p_t(x,y)\,dt\,\widetilde{\nu}(dx)\widetilde{\nu}(dy)\Big].$$
 Let us make the change of variables $u|x-y|^2=t$
 \begin{align}
\E^X\Big[&\int_{A\times A}\int_0^{|x-y|^2} \Ebr{\Bvec}{x}{y}{t}[e^{-\alpha F(t)}F(t)^{-s} ] p_t(x,y)\,dt\,\widetilde{\nu}(dx)\widetilde{\nu}(dy)\Big]\nonumber\\
=&\E^X\Big[\int_{A\times A}\int_0^{1} \Ebr{\Bvec}{x}{y}{u|x-y|^2}[e^{-\alpha F(u|x-y|^2)}F(u|x-y|^2)^{-s} ]  \frac{e^{-\frac{1}{2u}}}{2\pi u}\,du\,\widetilde{\nu}(dx)\widetilde{\nu}(dy)\Big]\nonumber\\
\leq& \int_{A\times A}\int_0^{1} \E^X\Big[\Ebr{\Bvec}{x}{y}{u|x-y|^2}[ \frac{e^{s\gamma X(x)-\frac{\gamma^2s^2}{2}\E^X[X^2(x)]}e^{s\gamma X(y)-\frac{\gamma^2s^2}{2}\E^X[X^2(x)]}}{F(u|x-y|^2)^{s}} ]  \Big]\frac{e^{-\frac{1}{2u}}}{2\pi u}\,du\, \nu (dx) \nu (dy).\label{petit}
\end{align}
We just have to analyze the expectation inside the integral. By translation invariance of the field $X$, it suffices to compute this expectation when $x=0$ and for a generic $y$. Thus we set
$$A(s,u,y)=  \E^X\Big[\Ebr{\Bvec}{0}{y}{u|y|^2}[ \frac{e^{s\gamma X(0)-\frac{\gamma^2s^2}{2}\E^X[X^2(0)]}e^{s\gamma X(y)-\frac{\gamma^2s^2}{2}\E^X[X^2(y)]}}{F(u|y|^2)^{s}} ]  \Big] .$$
By applying the same argument as in \cite[Lemma 28]{Rnew12} based on the Girsanov transform  and Kahane's convexity inequalities  
we may assume that the field $X$ is exactly scale invariant, meaning
that for all $\lambda\in [0,1]$ we have the following equality in law
\begin{equation}\label{exact}
(X(\lambda x))_{|x|\leq 4}  \stackrel{law}{=}(X(  x))_{|x|\leq 4} +\Omega_\lambda
\end{equation}
where $\Omega_\lambda$ is a centered Gaussian random variable independent of the whole field $(X(  x))_{|x|\leq 4}$ with variance $-\ln \lambda$. Let us denote by $K(x)=\E^X[X(x)X(0)]=\ln_+\frac{8}{|x|}$ its covariance kernel, which is nonnegative.

To use the scaling relation \eqref{exact}, we must force the Brownian bridge to stay within the ball $B(0,4)$. So we introduce the stopping time
$$S(W,t,r)=\inf\{ v\in [0,t]; W_v\not\in B(0,r)\}$$ and we have
 \begin{align}\label{porcback}
A(s,u,y)\le   \E^X\Big[\Ebr{\Bvec}{0}{y}{u|y|^2}[ \frac{e^{s\gamma X(0)-\frac{\gamma^2s^2}{2}\E^X[X^2(0)]}e^{s\gamma X(y)-\frac{\gamma^2s^2}{2}\E^X[X^2(y)]}}{F(u|y|^2\wedge S(\Bvec,u|y|^2,|y|))^{s}} ]  \Big] 
 \end{align}
Now we can use the  scaling relation \eqref{exact} together with that of the Brownian bridge to get by setting $e_y=y/|y|$. 
  \begin{align}\label{investA}
 &A(s,u,y)  \\
  &\leq u^{-s}\E\Big[\frac{e^{2s\gamma \Omega_{|y|}+\gamma^2s^2\ln|y|}}{|y|^{2s}e^{ s\gamma \Omega_{|y|}+\gamma^2s \ln|y|}}\Big]  \E^X\Big[ \Ebr{\Bvec}{0}{0}{1}[ \frac{e^{s\gamma X(0)-\frac{\gamma^2s^2}{2}\E^X[X^2(0)]}e^{s\gamma X(e_y)-\frac{\gamma^2s^2}{2}\E^X[X^2(e_y)]}}{\Big(\int_0^{1\wedge S(r\mapsto re_y+\sqrt{u}\Bvec_r,1,1)}e^{\gamma X(re_y+\sqrt{u} \Bvec_r )-\frac{\gamma^2}{2}\E^X[X^2]}\,dr\Big)^{s}} ]  \Big].\nonumber
 \end{align}
An elementary computation of the Laplace transform of Gaussian random variables gives the following explicit expression for the first expectation in the right-hand side of \eqref{investA}
\begin{equation}
\E\Big[\frac{e^{2s\gamma \Omega_{|y|}+\gamma^2s^2\ln|y|}}{|y|^{2s}e^{ s\gamma \Omega_{|y|}+\gamma^2s \ln|y|}}\Big]=|y|^{-\xi(s)}.
\end{equation}
The second expectation in the right-hand side of \eqref{investA} can be estimated with the Girsanov transform to get rid of the two exponentials in the numerator. Indeed we have the following elementary lemma, which is easily checked by computing the Laplace transform of the finite-dimensional distributions.

\begin{lemma}\label{L:girsanov}
Let $(X_t)_{t \in T}$ be a Gaussian process with covariance matrix $V = (v_{st})_{s,t \in T}$ on the probability space $(\Omega, \cF, \Pb)$.
\bi
\item
 Let $t_0$ be a point in $T$ and let $\tilde \Pb$ denote a probability measure with density $e^{\gamma X_{t_0} - \gamma^2 v_{t_0 t_0}/2 }$. Then the law of $X$ under $\tilde \Pb$ is that of $(X_t + \gamma v_{t_0t})_{t \in T}$. 
 \item Let now $t_0,t_1$ be two points in $T$ and let $\tilde \Pb$ be the probability measure with density
 \[
 e^{\gamma (X_{t_0}+X_{t_1}) - \gamma^2 (v_{t_0 t_0} +2 v_{t_0t_1}+ v_{t_1t_1} )/2}\,.
 \] Then the law of $X$ under $\tilde \Pb$ is that of $(X_t + \gamma (v_{t_0t}+v_{t_1t}))_{t \in T}$. 
 \ei
\end{lemma}

In order to use this Lemma in \eqref{investA}, one needs to add a factor $e^{\gamma^2 s^2 K(0,e_y)}$ which takes into account the covariance term. Note that this term is bounded by some constant $C$. 
Using this Lemma plus the non-negativity of the covariance kernel, we thus deduce 
   \begin{align*}
 A(s,u,y)& \\
  \leq &   \frac{C}{u^{s}|y|^{\xi(s)}}\E^X\Big[ \Ebr{\Bvec}{0}{0}{1}[ \Big(\int_0^{1\wedge S(r\mapsto re_y+\sqrt{u}\Bvec_r,1,1)}\!\!\!\!\!\!\!\!\!\!\!\!\!\!\!\!\!\!\!\!\!\!\!\!\!\!\!\!\!\!\!\!\!\!\!\!\!\!\!\!e^{\gamma X(re_y+\sqrt{u} \Bvec_r)-\frac{\gamma^2}{2}\E^X[X^2]+s\gamma^2 K(re_y+\sqrt{u} \Bvec_r)+s\gamma^2 K((1-r)e_y-\sqrt{u} \Bvec_r)}\,dr\Big)^{-s}]    \Big]\\
   \leq &    \frac{C}{u^{s}|y|^{\xi(s)}}\E^X\Big[ \Ebr{\Bvec}{0}{0}{1}[ \Big(\int_0^{1\wedge S(r\mapsto re_y+\sqrt{u}\Bvec_r,1,1)}e^{\gamma X(re_y+\sqrt{u} \Bvec_r)-\frac{\gamma^2}{2}\E^X[X^2]}\,dr\Big)^{-s}]    \Big].
\end{align*}

Let us now use the absolute continuity of Lemma \ref{absBB}. To this purpose, we first observe that the above expectation is smaller than the same one except that we replace the domain of integration over time $\int_0^{1\wedge S}$ by a smaller domain $\int_0^{(1/2)\wedge S}$. For this range of times $r\in[0,1/2\wedge S]$, we can apply Lemma \ref{absBB} and we see that the Radon-Nikodym derivative of $\Ebr{\Bvec}{0}{1}{1}$  with respect to $\E_0$ is bounded. 
Hence we get, 
after applying a change of variables $w=ur$,
\begin{align*}
A(s,u,y)
   \leq & C \frac{1}{u^{s}|y|^{\xi(s)}}\E_0 \Big[ \Big(\int_0^{(1/2) \wedge S(r\mapsto re_y+\sqrt{u}\Bvec_r,1,1)}e^{\gamma X(re_y+\sqrt{u} \Bvec_r)-\frac{\gamma^2}{2}\E^X[X^2]}\,dr\Big)^{-s}  \Big] \\
   =&C \frac{1}{ |y|^{\xi(s)}}  \E_0 \Big[  \Big(\int_0^{(u/2) \wedge S(w\mapsto wu^{-1}e_y+ \Bvec_w,1,1)}e^{\gamma X(wu^{-1}e_y+ \Bvec_w)-\frac{\gamma^2}{2}\E^X[X^2]}\,dw\Big)^{-s}  \Big] .
\end{align*}
In the above integral, the drift term $wu^{-1}e_y$ dominates the behaviour of the Brownian motion. To get rid of this effect, we will restrict the integration domain  of the integral inside the expectation to values that are very close to $0$. As we have $u\leq 1$, we deduce that $u^3/2\leq u/2$ in such a way that
$$A(s,u,y)\leq C |y|^{-\xi(s)} \E^X\big[ \Ebr{\Bvec}{0}{0}{1} \Big[ \Big(\int_0^{(u^3/2)\wedge S(w\mapsto wu^{-1}e_y+ B_w,1,u^{3/2})}e^{\gamma X(wu^{-1}e_y+ \Bvec_w)-\frac{\gamma^2}{2}\E^X[X^2]}\,dr\Big)^{-s}    \Big] \big].$$
Now we use the ordinary Cameron-Martin formula (Girsanov transform) for Brownian motion to get rid of the drift in the  latter quantity. This yields  
$$A(s,u,y)\leq C |y|^{-\xi(s)}\E_{0}\Big[ \Big(\int_0^{(u^3/2)\wedge S(w\mapsto  \Bvec_w,1,u^{3/2})}e^{\gamma X( \Bvec_w)-\frac{\gamma^2}{2}\E^X[X^2]}\,dr\Big)^{-s}e^{u^{-1}e_y\cdot \Bvec_{u^3}-\frac{1}{2}u}    \Big].$$
We can now easily get rid of each quantity in the above expectation. For instance, an abrupt use of the Cauchy--Schwarz inequality gives
\begin{align*}
A(s,u,y) \leq &C |y|^{-\xi(s)}\E_0\big[ F(u^3/2\wedge S(w\mapsto   \Bvec_w,1,u^{3/2}))^{-2s}     \big]^{1/2}e^{u/2}.
\end{align*}
By using the negative moments in \cite{GRV1} as well as their multifractal behaviour, we finally get (even if it means changing again the value of $C$)
\begin{equation}\label{estA}
A(s,u,y) \leq  C |y|^{-\xi(s)}u^{3\xi(-2s)/2}    \E_0\big[ F(1/2\wedge S(w\mapsto   \Bvec_w,1,1))^{-2s}     \big]^{1/2}\leq  C |y|^{-\xi(s)}u^{3\xi(-2s)/2}   
\end{equation}
so that plugging this estimate into \eqref{petit} yields
\begin{align*}
\E^X\Big[&\int_{A\times A}\int_0^{|x-y|^2} \Ebr{\Bvec}{x}{y}{t}[e^{-\alpha F(t)}F(t)^{-s} ] p_t(x,y)\,dt\,\widetilde{\nu}(dx)\widetilde{\nu}(dy)\Big]\\
\leq & C\int_{A\times A} \frac{1}{|x-y|^{\xi(s)}}\, \nu(dx) \nu (dy)\times \int_0^1u^{3\xi(-2s)/2-1} e^{-\frac{1}{2 u}}   \,du<+\infty.
\end{align*}

\subsubsection*{Long scales}

Now we focus on the long scale contribution
 $$B:=\E^X\Big[\int_{A\times A}\int_{|x-y|^{2 }}^\infty \Ebr{\Bvec}{x}{y}{t}[ e^{-\alpha F(t)}F(t)^{-s} ] p_t(x,y)\,dt\,\widetilde{\nu}(dx)\widetilde{\nu}(dy)\Big].$$
 For those $t$ that are larger than $|x-y|^2$, the Brownian bridge between $x$ and $y$ with lifetime $t$ behaves like a Brownian motion starting from $x$ between the times $0$ and $t/2$. Therefore it is not a bad idea to  apply Lemma \ref{absBB} to get rid of the bridge, and we get (recall that $\E_x = \E^X \otimes E^x_{\Bvec}$ is the annealed measure) 
 \begin{align}
B\leq & \E^X\Big[\int_{A\times A}\int_{|x-y|^{2}}^\infty \Ebr{\Bvec}{x}{y}{t}[ e^{-\alpha F(t/2)} F(t/2)^{-s} ] p_t(x,y)\,dt\,\widetilde{\nu}(dx)\widetilde{\nu}(dy)\Big]\label{largescale}\\
= & \E_x\Big[\int_{A\times A}\int_{|x-y|^{2}}^\infty e^{-\alpha F(t/2)} F(t/2)^{-s} 2 e^{\frac{|y-x|^2}{2t}-\frac{|\Bvec_{t/2}-x|^2}{t}} p_t(x,y)\,dt\,\widetilde{\nu}(dx)\widetilde{\nu}(dy)\Big]\nonumber\\
\leq  &C   \int_{A\times A}\int_{|x-y|^{2}}^\infty  \E_{x}\big[e^{-\alpha F(t/2)} F(t/2)^{-s}    e^{s\gamma X(x)+s\gamma X(y)- \gamma^2s^2 \E^X[X^2(x)]}\big] p_t(x,y)\,dt\, \nu(dx)\nu(dy). \nonumber
\end{align}
Now we want to separate in the above expectation the part that will serve to control the long time behaviour, i.e. the exponential weight $e^{-\alpha F (t/2)}$, and the part that will serve to control the scaling behaviour of the heat kernel, i.e. $F(t/2)^{-s}$. To this purpose, we first use the usual trick of replacing the MFF $X$ by the exactly scale invariant kernel with decorrelation cutoff at length $2$, still denoted by $X$ (we can do this as in  \cite[Lemma 28]{Rnew12} because the mapping $x\mapsto e^{-\alpha x}x^{-s}$ is convex). Its covariance kernel is given by $K(x,y)=\ln_+\frac{2}{|x-y|}$ and has the scaling relation \eqref{exact} for $|x|\leq 1$. Then we introduce the functionals
\begin{equation}
F^c_r(t)=\int_0^t\ind_{\{\Bvec_v\not \in B(0,r)\}}\,F(dv),\quad \quad F^i(t)=\int_0^t\ind_{\{\Bvec_v  \in B(0,1)\}}\,F(dv).
\end{equation}
Notice that $F^c_r(t)+F^i(t)\leq F(t)$ for all $t$ and $r\geq 1$. Finally we set
$$B(t,x,y):=  \E_{x}\Big[e^{-\alpha F (t/2)} F(t/2)^{-s}   e^{s\gamma X(x)+s\gamma X(y)- \gamma^2s^2 \E^X[X^2(x)]}\Big] $$
Since $F(t/2)\geq F(|y-x|^2/2)$ for $t\geq |x-y|^2$, we have
 \begin{align*}
B(t,x,y)\leq &    \E_{x}\Big[e^{-\alpha (F^c_2(t/2)-F^c_2(|x-y|^2/2))} F^i(|x-y|^2/2)^{-s}    e^{s\gamma X(x)+s\gamma X(y)- \gamma^2s^2 \E^X[X^2(x)]}\Big] .
\end{align*}
Notice that under the annealed probability measure $  \Pb_x$ 
and for $x,y\in A\subset B(0,1/2)$, the random variables $F^c_2(t/2)-F^c_2(|x-y|^2/2)$ and $F^i(|x-y|^2/2)^{-s}   e^{s\gamma X(x)+s\gamma X(y)- \gamma^2s^2 \E^X[X^2(x)]}$ are independent (use the decorrelation cutoff of the field $X$ and the strong Markov property of the Brownian motion) so that we get
 \begin{align*}
B(t,x,y)\leq &     \E_x\Big[e^{-\alpha ( F^c_2(t/2)-F^c_2(|x-y|^2/2))}\Big]  \E_x\Big[ F^i(|x-y|^2/2)^{-s}    e^{s\gamma X(x)+s\gamma X(y)- \gamma^2s^2 \E^X[X^2(x)]}\Big] .
\end{align*}
Now we use \eqref{estA} with $u=1/2$ to estimate the second expectation and get (actually, one has to re-derive  \eqref{estA} with $F$ replaced by $F^i$ but this is harmless)
\begin{equation}\label{estB}
B(t,x,y)\leq    C   \E_x\Big[e^{-\alpha  ( F^c_2(t/2)-F^c_2(|x-y|^2/2))}\Big] \frac{1}{|x-y|^{\xi(s)}}.
\end{equation}
It remains to treat the expectation $   \E_x[e^{-\alpha  ( F^c(t/2)-F^c_2(|x-y|^2/2))}]  $ and our purpose is to establish that
\begin{equation}\label{estexp}
  \E_x[e^{-\alpha  ( F^c_2(t/2)-F^c_2(|x-y|^2/2))}]   \leq C(t-|x-y|^2)^{-\delta}
\end{equation}
uniformly with respect to $x,y\in A$ and where $\delta>0$ is an in \eqref{eucap0}. By using the Markov property of the Brownian motion at time $|x-y|^2$, the stationarity of the field $X$ and the fact that $B(z,2)\subset B(0,4)$ for all $z$ such that $|z|\le 2$ we deduce
\begin{equation}\label{estexp2}
\E_x\big[e^{-\alpha  ( F^c_2(t/2)-F^c_2(|x-y|^2/2))}\big] \leq \E_0\big[e^{-\alpha  F^c_4(t/2-|x-y|^2/2)}\big].
\end{equation}
Let us introduce the stopping times $T_0=0$  and for $n\geq 0$
$$T_{n+1}=\inf\{v>\bar{T}_{n};|\Bvec_v|=4n+10\}\quad \quad  \bar{T}_n=\inf\{v>T_{n};|\Bvec_v-\Bvec_{T_n}|=1\}.$$
Under $\Pb_0$ the random variables  $(F^c_4(\bar{T}_n)- F^c_4(T_n))_n  $ are independent and identically distributed. Furthermore $F^c_4(\bar{T}_n)- F^c_4(T_n)=F(\bar{T}_n)- F^c_4(T_n)$. We deduce
 \begin{align*}
 \E_0[e^{-\alpha  F^c_4(t/2-|x-y|^2/2)}]\leq & \E_0[e^{-\alpha  \sum_{n,\bar{T}_n\leq (t-|x-y|^2)/2} (F(\bar{T}_n)- F^c_4(T_n))}]\\
 \leq & \E_0[e^{-\alpha  \sum_{n=1}^K (F(\bar{T}_n)- F^c_4(T_n))}]+\Pb_0(\sup_{v\leq (t-|x-y|^2)/2}|\Bvec_v|<4K+10)\\
 \leq &  \E_0[e^{-\alpha   F(\bar{T}_0)}]^K+C\frac{K^{1/2}}{(t-|x-y|^2)^{1/4} }
 \end{align*}
In the last line, we have used standard estimates of the exit times of the Brownian motion out of the balls  in a way that is quite non optimal but this is enough for our purposes. It suffices to take $K=(t-|x-y|^2)^{1/2 - 2\delta }$ and to observe that $\E_0[e^{-\alpha   F(\bar{T}_0)}]<1$, and we obtain \eqref{estexp}.
Now, by making the change of variables $u=\frac{t}{|x-y|^2}$ , we obtain
\begin{align*}
\int_{|x-y|^{2}}^\infty  (t-|x-y|^2)^{-\delta} p_t(x,y) \,dt& = \int_{|x-y|^{2}}^\infty  (t-|x-y|^2)^{-\delta} e^{-\frac{|x-y|^2}{2t}} \,\frac{dt}{t}\\
&=  |x-y|^{-2\delta}\int_{1}^\infty  (1-\frac{1}{u})^{-\delta} e^{-\frac{1}{2u}} \,\frac{du}{u^{1+\delta}}\\
&= C   |x-y|^{-2\delta}.
\end{align*}

By gathering  \eqref{estB}+\eqref{estexp} and plugging the result into \eqref{largescale} we obtain
 \begin{align*}
 B\leq & C \int_{A\times A}\frac{1}{|x-y|^{\xi(s)}}\Big(\int_{|x-y|^{2}}^\infty  (t-|x-y|^2)^{-\delta} p_t(x,y)\,dt\Big)\, \nu(dx)\nu(dy)\\
 \leq & C \int_{A\times A}\frac{}{|x-y|^{\xi(s) + 2\delta }} \nu(dx)\nu(dy) <\infty\,. 
 \end{align*}
 This completes the proof.\qed  
\appendix

\section{(No-)doubling property of the Liouville measures}\label{sec:doubling}
In this appendix, we shall prove the following Theorem which illustrates that Liouville measures possess very poor doubling properties:

\begin{theorem}\label{doubling}
Let $\gamma<\gamma_c=2$ and consider a MFF $X$ in the plane. Set
$$\eta_c=\frac{\gamma^2}{4+\gamma^2}.$$
For any $ \eta\in]\eta_c,1]$, there is a constant $C=C(X,\gamma,\eta)$ which is a.s. finite such that for any $x\in[0, 1]^2$ and any radius $r\in(0,1]$, one has 
\[
M_\gamma(B(x,2r)) \leq C\, M_\gamma(B(x,r))^{1-\eta}\,.
\]
\end{theorem}

\begin{remark}\label{}
It can be seen that $\eta_c=\eta_c(\gamma)$ is optimal in the sense that 
$$\underset{r \to 0}{\overline{\lim}}\sup_{x\in B(0,1)}\frac{M_\gamma(B(x,2r))}{M_\gamma(B(x,r))^{1-\eta}}=\left\{\begin{array}{ll}
<+\infty, & \text{ if }\eta>\frac{\gamma^2}{4+\gamma^2}\\
=+\infty, & \text{ if }\eta<\frac{\gamma^2}{4+\gamma^2}.
\end{array}\right.$$
We provide a rigorous proof of the sharpness of $\eta_c$ in this case of dyadic cascades. See Theorem \ref{doubling2} below. 
\end{remark}

\ni
{\it Proof}. Let $(X_\epsilon)_\epsilon$ be a white noise decomposition of the massive free field, meaning that $(X_\epsilon)_\epsilon$ is a family of centered Gaussian fields with covariance given by 
\begin{align}\label{e.cov}
\E[X_\epsilon(x)X_{\epsilon'}(y)]=\int_1^{\frac{1}{\epsilon\wedge \epsilon '}}\frac{k(u(x-y))}{u}\,du
\end{align}
for $\epsilon,\epsilon'\in]0,1]$ and  $k(z)=\frac{1}{2}\int_0^\infty e^{-\frac{m^2}{2v}|z|^2-\frac{v}{2}}\,dv$. Notice that we skip the index $X$ in the expectation as we only consider expectations with respect to the family $(X_\epsilon)_\epsilon$ in the following (no distinction with the Brownian motion  is needed).

We want to evaluate the quantity   
\begin{align}\label{doub1}
\Pb\Big(\frac{M_\gamma(B(0,16r))}{M_\gamma (B(0,r))^{1-\eta}}\geq 1\Big) .
\end{align}
Recall that the massive free field is star scale invariant in the sense of \cite{Rnew1}, meaning that for $\epsilon>0$ we have the following equality in law
$$M_\gamma(dz)=e^{\gamma X_\epsilon(z)-\frac{\gamma^2}{2}\E[X_\epsilon(z)^2]}\epsilon^2M^\epsilon_\gamma(dz/\epsilon)$$
where $M^\epsilon_\gamma$ is a random measure independent of $X_\epsilon$ with the same law as $M_\gamma$. So we can write
\begin{align*}
 M_\gamma(B(0,r)) =&r^2e^{\gamma X_r(0)-\frac{\gamma^2}{2}\E[X_r(0)^2]}\int_{B(0,r)}e^{\gamma Y_r(z)}M^r_\gamma(dz)
 \end{align*}
 where $Y_r(z)=X_r(z)-X_r(0)$ is independent of $M^r_\gamma$. Therefore  
 $$M_\gamma(B(0,r))\geq r^2e^{\gamma X_r(0)-\frac{\gamma^2}{2}\E[X_r(0)^2]}e^{\gamma  \min_{z\in B(0,r)}Y_r(z)}M^r_\gamma(B(0,1)).$$ In the same way, 
 $$M_\gamma(B(0,16r))\leq r^2e^{\gamma X_r(0)-\frac{\gamma^2}{2}\E[X_r(0)^2]}e^{\gamma \max_{z\in B(0,16r)}Y_r(z)}M^r_\gamma(B(0,16)).$$
We deduce for some $ \delta_1,\delta_2>0$ (to be determined later)
\begin{align}
\Pb\Big(&\frac{M_\gamma(B(0,16r))}{M_\gamma (B(0,r))^{1-\eta}} \geq 1\Big)\nonumber\\
\leq & \Pb\Big(r^{(2+\frac{\gamma^2}{2})\eta}e^{\eta\gamma X_r(0)}e^{\gamma \max_{z\in B(0,16r)} Y_r(z)}M^r_\gamma(B(0,16))\geq  e^{(1-\eta)\gamma\min_{z\in B(0,r)}Y_r(z)}M^r_\gamma(B(0,1))^{1-\eta} \Big) \nonumber\\
\leq &  \Pb\Big(r^{(2+\frac{\gamma^2}{2})\eta}e^{\eta\gamma X_r(0)}r^{-\gamma \delta_1}M^r_\gamma(B(0,16))\geq  r^{(1-\eta)\gamma \delta_2}M^r_\gamma(B(0,1))^{1-\eta} \Big)  \nonumber\\
&+\Pb\Big(\max_{z\in B(0,16r)} Y_r(z) \geq -\delta_1\ln r  \Big)+\Pb\Big(\min_{z\in B(0,16r)} Y_r(z) \leq \delta_2\ln r  \Big)
\label{doub2}.
\end{align}
Let us denote by $F$ the function defined for    $x\in\R$ by
$$F(x)=\frac{1}{\sqrt{2\pi}}\int_x^{\infty} e^{-\frac{u^2}{2}}\,du$$
and by $Z$ the random variable $Z=\frac{M_\gamma(B(0,2))}{M_\gamma(B(0,1))^{1-\eta}} $. 
We have for any $a>0$
\begin{align}
 \Pb\Big(&r^{(2+\frac{\gamma^2}{2})\eta}e^{\eta\gamma X_r(0)}r^{-\gamma \delta_1}M^r_\gamma(B(0,16))\geq  r^{(1-\eta)\gamma \delta_2}M^r_\gamma(B(0,1))^{1-\eta} \Big)  \nonumber\\
\leq &  \Pb\Big(  X_r(0)  \geq  (\eta\gamma)^{-1}\big( (2+\frac{\gamma^2}{2}) \eta   - (1-\eta)\gamma \delta_2 - \gamma\delta_1\big)\ln\frac{1}{r} + (\eta\gamma)^{-1}\ln Z^{-1} \Big)  \nonumber \\
=& F\Big((\eta\gamma)^{-1}\big((2+\frac{\gamma^2}{2})\eta-(1-\eta)\gamma \delta_2-\gamma\delta_1\big)(-\ln r)^{1/2} + (\eta\gamma)^{-1}(-\ln r)^{-1/2}\ln Z^{-1} \Big)  \nonumber \\
\leq& F\Big((\eta\gamma)^{-1}\big((2+\frac{\gamma^2}{2})\eta-(1-\eta)\gamma \delta_2-\gamma\delta_1-a\big)(-\ln r)^{1/2}  \Big)\nonumber\\
&+\Pb(\ln Z^{-1}\leq a \ln r).\label{doub3}
\end{align}
Let us estimate the last probability in the latter expression. For $q<4/\gamma^2$ and $p,\bar{p}>1$ such that $ p^{-1}+\bar{p}^{-1}=1$ and $pq<4/\gamma^2$, we have by using in turn the Markov inequality and the H\"older inequality
\begin{align}
\Pb(\ln Z^{-1}\leq a \ln r)&=\Pb( Z\geq  r^{-a})\nonumber\\
&\leq r^{aq}\E[Z^q]\nonumber\\
&\leq r^{aq}\E[M_\gamma(B(0,16))^{pq}]^{1/p}\E[M_\gamma(B(0,1))^{-\bar{p}q(1-\eta)}]^{1/\bar{p}}.\label{doub4}
\end{align}
Because of our choice of $p,q$ the first expectation $\E[M_\gamma(B(0,16))^{pq}]$ is finite (existence of moments up to $4/\gamma^2$ \cite{cf:Kah,review}) as well as the second $\E[M_\gamma(B(0,1))^{-\bar{p}q(1-\eta)}]$ (existence of negative moments of all orders \cite{review}).

Finally we estimate the probabilities involving the $\max$ and $\min$ in \eqref{doub2}. The key point is to observe that the Gaussian process $Y_r$ does not fluctuate too much in such a way that its maximum (and minimum) possesses a   Gaussian right tail distribution.   From the covariance structure ~\eqref{e.cov}, it is easy to observe that
$$\max_{z\in B(0,16r)}\gamma Y_r(z) =\max_{z\in B(0,16)}\gamma Y_r(r z) \quad\text{ and }\quad \forall z,z'\in B(0,16),\quad \E[(Y_r(r z) -Y_r(r z') )^2]\leq C|z-z'|.$$
Using for example \cite{Ledoux}(Thm. 7.1, Eq. (7.4)), one can then deduce  
$$\forall x\geq 1,\quad \sup_r\Pb(\max_{z\in B(0,16r)} Y_r(z)\geq x)\leq Ce^{-cx^2}$$ for some constants $C,c>0$ which does not depend on $r$ (a similar relation holds for the left tail distribution of the minimum).

Now we gather the previous relations and to fix the values of the parameters $\delta,\delta_1,\delta_2,a,q$. Let us choose $a'$ such that
$$0<a'<(2+\gamma^2/2)\eta\quad \text{and} \quad a'\frac{4}{\gamma^2}>2.$$
(This is possible since $\eta > \eta_c$.) By continuity, we can find $ \delta_1,\delta_2,a,q>0$ such that
$$(2+\frac{\gamma^2}{2})\eta-(1-\eta)\gamma \delta_2-\gamma\delta_1-a>0\quad \text{and} \quad aq>2. $$
By gathering  \eqref{doub2}+\eqref{doub3}+\eqref{doub4} and standard estimates on the asymptotic behaviour of $F$ for large $x$, we get for some $c>0$
\begin{align*}
\Pb\Big(&\frac{M_\gamma(B(0,16r))}{M_\gamma (B(0,r))^{1-\eta}} \geq 1\Big)\leq ce^{-\frac{(2+\frac{\gamma^2}{2} -(1-\eta)\gamma \delta_2-\gamma\delta_1-a)^2}{2\gamma^2\eta^2}(\ln r)^2}+cr^{aq}+e^{-c\delta_1^2 (\ln r)^2}+e^{-c\delta_2^2 (\ln r)^2}
\end{align*}
in such a way that that for some $\delta'>0$
\begin{equation}\label{doubf}
\sup_{r\in]0,1]} r^{-2-\delta'}\Pb\Big(\frac{M_\gamma(B(0,16r))}{M_\gamma (B(0,r))^{1-\eta}} \geq 1\Big)<\infty.
\end{equation}

It remains to conclude. Let us consider the square $S=[0,1]^2$ and let us denote by $D_n$ the dyadic numbers of the square of order $n$, i.e. of the form  $x=(\frac{k}{2^n},\frac{k'}{2^n})$ for some $k,k'\in [0,2^n]\cap\N$. We have by invariance under translations of $M_\gamma$ and \eqref{doubf}
\begin{align*}
\Pb\Big(\max_{x\in D_n}\frac{M_\gamma(B(x,16\times 2^{-n}))}{M_\gamma (B(x,2^{-n}))^{1-\eta}} \geq 1\Big)\leq &\sum_{x\in D_n}\Pb\Big( \frac{M_\gamma(B(x,16\times 2^{-n}))}{M_\gamma (B(x,2^{-n}))^{1-\eta}} \geq 1\Big)\\
=&2^{2n}\Pb\Big( \frac{M_\gamma(B(0,16\times 2^{-n}))}{M_\gamma (B(0,2^{-n}))^{1-\eta}} \geq 1\Big)\\
\leq & 2^{-n\delta'}.
\end{align*}
By using the Borel-Cantelli lemma, we deduce that there exists a random constant $C$ such that
$$\sup_n\max_{x\in D_n}\frac{M_\gamma(B(x,16\times 2^{-n}))}{M_\gamma (B(x,2^{-n}))^{1-\eta}}\leq C.$$
Finally we consider any $y\in S$ and $r>0$. Let us consider $n$ such that $2^{-(n+1)}\leq r <2^{-n}$. Let us denote by $x$ the point in $D_{n+2}$ which is the closest to $y$ (in particular, notice that $|x-y|\leq 2^{-(n+2)}$). We have
\begin{align*}
M_\gamma(B(y,2 r))\leq & M_\gamma(B(x,4\times 2^{-n}))\\
\leq & C M_\gamma(B(x,2^{-(n+2)}))^{1-\eta}\\
\leq & C M_\gamma(B(y,r))^{1-\eta},
\end{align*}
which completes the proof of the second relation. \qed

\medskip
We will prove now that the threshold given in Theorem \ref{doubling} is sharp in the context of lognormal multiplicative cascades (see \cite{KP} for further details).  Let us recall briefly the setup. If $x\in[0,1]$, we will denote by $(x_k)_{k\geq 1}\in\{0;1\}^{\N^*}$ its dyadic decomposition, i.e.
 $$x=\sum_{k=1}^\infty\frac{x_k}{2^k}.$$
$D_n\subset [0,1] $ stands for the set of dyadic numbers of order $n$ $D_n=\big\{\frac{k}{2^n};k=0,\dots,2^n\big\}$. We set $D=\bigcup_n D_n $. We denote by  $\pi_n:[0,1]\to D_n$  the projection of $[0,1] $ onto $D_n$, i.e.
 $$\forall x\in [0,1],\quad \pi_n(x)= \sum_{k=1}^n\frac{x_k}{2^k}.$$ 
Let us set $I_n^k=[\frac{k}{2^n},\frac{k+1}{2^n}[$ for $k=0,\dots,2^{n}-1$ and for $x\in [0,1[$, we denote by $I_n(x)$ the unique interval among the family $(I_n^k)_{k=0,\dots,2^{n}-1}$ containing $x$.

Consider  a sequence of i.i.d.  random variables $(X_i)_{i\in D}$ indexed by the set $D=\bigcup_nD_n$ of dyadic numbers of $[0,1]$ with common law that of a standard Gaussian random variable.

Finally, for $\gamma>0$ such that $\gamma^2 < 2 \ln 2$, we consider the multiplicative cascade measure $\mu$  defined by 
\begin{equation*}
\mu(dx) = \underset{n \to \infty}{\lim}   e^{\gamma X_{\pi_1(x)}+ \cdots +\gamma X_{\pi_n(x)} -\gamma^2/2 n }\,dx.
\end{equation*}
For $\eta>0$, we want to study the maximum of $\frac{\mu(I_n(y))}{\mu(I_{n+1}(y))^{1-\eta}}$ for $y\in [0,1]$. 
\begin{theorem}\label{doubling2}
If $\eta<\frac{\gamma^2}{\gamma^2+2\ln 2}$, then for all $R>0$
$$\lim_{n\to\infty}\Pb\Big(\max_{y\in[0,1]}\frac{\mu(I_n(y))}{\mu(I_{n+1}(y))^{1-\eta}}\geq R\Big)=1.$$
\end{theorem}

\noindent {\it Proof.} In the context of multiplicative cascades,  there exists an independent family of random variables $(Z_{y}^{(n+1)})_{y\in D_{n+1}}$   with law $\mu[0,1]$ and such that  for all   $x\in [0,1]$
\begin{equation*}
\mu(I^{n+1}(x))= \frac{1}{2^{n+1}}e^{\gamma X_{\pi_1(x)}+ \cdots +\gamma X_{\pi_{n+1}(x)} -\gamma^2/2 (n+1) } Z_{\pi_{n+1}(x)}^{(n+1)}.
\end{equation*}   
The total mass $\mu[0,1]$ is such that there exists $\alpha>0$, more precisely $\alpha=\frac{2 \ln 2}{\gamma^2}$, such that $\Pb( \mu [0,1]>u ) \sim \frac{1}{u^\alpha}$ for $u\to\infty$. Hence, from the Borel-Cantelli lemma, we get that for all $\delta>0$, there exists a random constant $c$ such that
\begin{equation}\label{limmax}
\lim_{n\to\infty}\Pb\big(\max_{ y\in D_{n+1} }  Z_{y}^{(n+1)} \geq   2^{n/\alpha-\delta}\big)=1  .
\end{equation}
For $y\in D_{n+1}$ we denote by $y^c$ the unique element of $D_{n+1}$ such that $\pi_n(y)=\pi_n(y^c)$ and $y\not =y^c$.  Furthermore, for $y\in [0,1[$, we set $Z_{n}(y)=  \frac{1}{2^{n}}e^{\gamma X_{\pi_1(y)}+ \cdots +\gamma X_{\pi_{n}(y)} -\gamma^2/2 n } $. 
Then we have for all $R>0$
\begin{align*}
\Pb\Big(&\max_{y\in D_{n+1}}\frac{\mu(I_n(y))}{\mu(I_{n+1}(y^c))^{1-\eta}}\geq R\Big)\\
=& \Pb\Big(\max_{y\in D_{n+1}}2^\eta Z_{n}(y)^\eta  \frac{e^{  \gamma X_{\pi_{n+1}(y) }-\gamma^2/2} Z_{\pi_{n+1}(y) }^{(n+1)}    + e^{\gamma X_{{\pi_{n+1}(y^c) }} -\gamma^2/2} Z_{\pi_{n+1}(y^c) }^{(n+1)}     }{  (Z_{\pi_{n+1}(y^c) }^{(n+1)})^{1-\eta} }\geq R\Big) \\
\geq & \Pb\Big(\max_{y\in D_{n+1}}  Z_{n}(y)^\eta  \frac{e^{  \gamma X_{\pi_{n+1}(y) }-\gamma^2/2} 2^{n/\alpha-\delta}     }{  (Z_{\pi_{n+1}(y^c) }^{(n+1)})^{1-\eta} }\geq 2^{-\eta}R2^{-n/\alpha+\delta} \Big)- \Pb\Big(\max_{ y\in D_{n+1} }  Z_{y}^{(n+1)}  \leq 2^{n/\alpha-\delta} \Big) .
\end{align*}
A simple computation on the maximum of Gaussian random variables (similar to the proof of Theorem \ref{doubling}) shows that for $\eta$ such that $\eta< \frac{1}{\alpha (  1+\gamma^2/(2 \ln 2) )}=\frac{\gamma^2}{\gamma^2+2\ln 2}$, we have
$$\lim_{n\to\infty} \Pb\Big(\max_{y\in D_{n+1}}  Z_{n}(y)^\eta  \frac{e^{  \gamma X_{\pi_{n+1}(y) }-\gamma^2/2} 2^{n/\alpha-\delta}     }{  (Z_{\pi_{n+1}(y^c) }^{(n+1)})^{1-\eta} }\geq 2^{-\eta}R2^{-n/\alpha+\delta} \Big)=1.$$
We conclude with \eqref{limmax}.\qed

 \end{document}